\documentclass[12pt,leqno]{article}
\usepackage{amssymb,graphics,latexsym}
\usepackage{graphicx}
\usepackage{amsmath}
\usepackage{amsthm}
\large\normalsize

\def\R{\hbox{\bf\rlap{I}{\hbox to 2 pt{}}R}}
\def\tb#1#2{\mathop{#1\vphantom{\sum}}\limits_{\displaystyle #2}}

\newcommand{\re}{{\rm Re\, }}
\newcommand{\im}{{\rm Im\, }}
\newcommand{\tr}{{\rm tr\, }}
\newcommand{\dia}{{\rm diag\, }}

\newcommand{\rank}{{\rm rank\, }}

\begin{document}
\thispagestyle{empty}
\begin{center}
\section*{Numerical Ranges of KMS Matrices}

\vspace*{10mm}

{\bf Hwa-Long Gau$^*$}\hspace{.5cm} and   \hspace{.5cm}{\bf Pei Yuan Wu$^{**}$}

\vspace{10mm}

In memory of B\'{e}la Sz\H{o}kefalvi-Nagy on his 100th anniversary

\end{center}

\vspace{10mm}

\noindent
{\bf Abstract}

A KMS matrix is one of the form
$$J_n(a)=\left[\begin{array}{ccccc} 0 & a & a^2 & \cdots & a^{n-1}\\ & 0 & a & \ddots & \vdots\\ & & \ddots & \ddots & a^2\\ & & & \ddots & a\\ 0 & & & & 0\end{array}\right]$$
for $n\ge 1$ and $a$ in $\mathbb{C}$. Among other things, we prove the following properties of its numerical range: (1) $W(J_n(a))$ is a circular disc if and only if $n=2$ and $a\neq 0$, (2) its boundary $\partial W(J_n(a))$ contains a line segment if and only if $n\ge 3$ and $|a|=1$, and (3) the intersection of the boundaries $\partial W(J_n(a))$ and $\partial W(J_n(a)[j])$ is either the singleton $\{\min\sigma(\re J_n(a))\}$ if $n$ is odd, $j=(n+1)/2$ and $|a|>1$, or the empty set $\emptyset$ if otherwise, where, for any $n$-by-$n$ matrix $A$, $A[j]$ denotes its $j$th principal submatrix obtained by deleting its $j$th row and $j$th column ($1\le j\le n$), $\re A$ its real part $(A+A^*)/2$, and $\sigma(A)$ its spectrum.

\vspace{1cm}



\vspace{.5cm}
\noindent
${}^*$Partially supported by the National Science Council of the
Republic of China under project NSC 101-2115-M-008-006.

\noindent
${}^{**}$Partially supported by the National Science Council of the
Republic of China under project NSC 101-2115-M-009-004 and by the MOE-ATU.

\newpage

\noindent
{\bf\large 1. Introduction}

\vspace{.3cm}

An $n$-by-$n$ matrix of the form
$$J_n(a)=\left[\begin{array}{ccccc} 0 & a & a^2 & \cdots & a^{n-1}\\ & 0 & a & \ddots & \vdots\\ & & \ddots & \ddots & a^2\\ & & & \ddots & a\\ 0 & & & & 0\end{array}\right],$$
where $n\ge 1$ and $a$ is in $\mathbb{C}$, is called an (upper-triangular) \emph{KMS matrix}. Although in a very simple form, such matrices, as a meeting ground of nilpotent matrices, Toeplitz matrices, nonnegative matrices, $S_n$-matrices and $S_n^{-1}$-matrices, have many interesting properties, especially concerning their numerical ranges. Recall that the \emph{numerical range} $W(A)$ of an $n$-by-$n$ matrix $A$  is the set $\{\langle Ax, x\rangle : x\in \mathbb{C}^n, \|x\|=1\}$, where $\langle\cdot,\cdot\rangle$ denotes the standard inner product in $\mathbb{C}^n$ and $\|\cdot\|$ its associated norm. It is known that $W(A)$ is a compact convex subset of the complex plane. For other properties of the numerical range, the reader may consult \cite[Chapter 1]{14}.

\vspace{.5cm}

The study of the numerical range of $J_n(a)$ was started by Gaaya in \cite{2, 3}. If $\re A$ denotes the real part $(A+A^*)/2$ of a matrix $A$ and $I_n$ the $n$-by-$n$ identity matrix, then, for $0\le a<1$,
$$2\re J_n(a)+I_n=\left[\begin{array}{cccc}  1 & a & \cdots & a^{n-1}\\ a & 1 & \ddots & \vdots\\ \vdots & \ddots & \ddots & a\\ a^{n-1} & \cdots & a & 1\end{array}\right]$$
is a Toeplitz matrix associated to the Poisson kernel $P_a(e^{it})=(1-a^2)/|1-ae^{it}|^2$, first introduced by Kac, Murdock and Szeg\H{o} \cite{15}, with eigenvalues $P_a(e^{it_k^{(n)}})$, where $t_k^{(n)}$, $1\le k\le n$, are the roots of the equation $\sin((n+1)t)-2a\sin(nt)+a^2\sin((n-1)t)=0$ (cf. \cite[pp. 69--70]{12}). This was exploited in \cite{2,3} to obtain the value of  the numerical radius of $J_n(a)$. The purpose of this paper is to launch a systematic study of properties of the numerical range of $J_n(a)$. More specifically, we determine, in Section 2 below, when $W(J_n(a))$ is a circular disc, namely, we prove that this is the case if and only if $n=2$ and $a\neq 0$ (Theorem 2.3). Two proofs of this will be given. One involves computations with the Kippenhahn polynomial of $J_n(a)$ while the other relates to the theories of $S_n$-, $S_n^{-1}$- and nonnegative matrices. Both of them make use of B\'{e}zout's theorem from algebraic geometry. The relations to matrices of classes $S_n$ and $S_n^{-1}$ also help us in showing that the boundary of $W(J_n(a))$ contains a line segment if and only if $n\ge 3$ and $|a|=1$ (Proposition 2.11), thus generalizing the corresponding result for $a=1$ (cf. \cite[Lemma 1.3 (4)]{10}). The main concern in Section 3 is the relations between the numerical ranges of $J_n(a)$ and its compressions. Recall that an $m$-by-$m$ matrix $A$ is a \emph{compression} of an $n$-by-$n$ matrix $B$ if $A=V^*BV$ for some $n$-by-$m$ matrix $V$ with $V^*V=I_m$. Equivalently, this is the same as requiring that $B$ be unitarily similar to a matrix of the form {\scriptsize$\left[  \begin{array}{cc}   A & \ast \\ \ast & \ast \\ \end{array}  \right]$}. Obviously, in this case we have $W(A)\subseteq W(B)$. We will prove that if $1\le m<n$ and $a\neq 0$, then every $m$-by-$m$ compression $A$ of $J_n(a)$ is such that $W(A)\subsetneqq W(J_n(a))$ (Proposition 3.1), and, for the special case of $|a|\neq 0, 1$ and $J_n(a)$ unitarily similar to {\scriptsize$\left[  \begin{array}{cc}   A & \ast \\ 0 & \ast \\ \end{array}  \right]$}, we even have the containment of $W(A)$ in the interior of $W(J_n(a))$ (Proposition 3.5). Finally, for an $(n-1)$-by-$(n-1)$ principal submatrix $A$ of $J_n(a)$, we determine the intersection $\partial W(A)\cap\partial W(J_n(a))$, which turns out to consist of at most one element (Theorem 3.8).

\vspace{5mm}

We conclude this section with some notations frequently used in the discussions below. For an $n$-by-$n$ matrix $A$, we use $A^T$, $\re A$, $\im A$, $\tr A$, $\det A$ and $\rank A$ to denote its transpose, real part $(A+A^*)/2$, imaginary part $(A-A^*)/(2i)$, trace, determinant and rank, respectively. The spectrum and numerical radius of $A$ are $\sigma(A)$ and $w(A)$ ($\equiv\max\{|z|: z\in W(A)\}$), respectively. The $n$-by-$n$ diagonal matrix with diagonal entries $a_1, \ldots, a_n$ is denoted by $\dia(a_1, \ldots, a_n)$, and the $n$-by-$n$ zero matrix (respectively, identity matrix) is $0_n$ (respectively, $I_n$). The subspace generated by the vectors in $\mathcal{S}\subseteq \mathbb{C}^n$ (or the span of $\mathcal{S}$) is $\bigvee\mathcal{S}$. For a subset $\bigtriangleup$ of $\mathbb{C}$, $\bigtriangleup^{\wedge}$ (respectively, $\#\bigtriangleup$) denotes the convex hull (respectively, cardinal number) of $\bigtriangleup$.

\vspace{1cm}

\noindent
{\bf\large 2. Circular Disc and Line Segment}

\vspace{.3cm}

We start with the following proposition, which gives some basic properties of $J_n(a)$ and its numerical range.

\vspace{.5cm}

{\bf Proposition 2.1.} (a) \emph{If} $|a|=|b|$, \emph{then $J_n(a)$ and $J_n(b)$ are unitarily similar}.

(b) \emph{If $a$ and $b$ are nonzero}, \emph{then $J_n(a)$ and $J_n(b)$ are similar}.

(c) \emph{$W(J_n(a))$ is symmetric with respect to the $x$-axis}.

(d) \emph{For any $n\ge 2$ and $a\neq 0$}, \emph{$J_n(a)$ is irreducible}, \emph{$0$ is in the interior of $W(J_n(a))$ and $\partial W(J_n(a))$ is a differentiable curve}.

(e) \emph{If $|a|\le|b|$}, \emph{then $W(J_n(a))$ is contained in $W(J_n(b))$}.

\vspace{.5cm}

Recall that a square matrix is said to be \emph{irreducible} if it is not unitarily similar to the direct sum of two other matrices. To prove (e) of the preceding proposition, we need the next lemma, which is a generalization of \cite[Lemma 3.2]{11}.

\vspace{.5cm}

{\bf Lemma 2.2.} \emph{If $|a_j|\le |b_j|$ for $1\le j\le n-1$}, \emph{and} $A_{ij}$, $1\le i\le j\le n$, \emph{are operators} (\emph{on appropriate spaces}), \emph{then}
$$W(\left[\begin{array}{cccc} A_{11} & a_1A_{12} & \cdots & a_1\cdots a_{n-1}A_{1n}\\ & A_{22} & \ddots & \vdots\\ & & \ddots & a_{n-1}A_{n-1, \, n}\\ 0 & & & A_{nn}\end{array}\right])\subseteq W(\left[\begin{array}{cccc} A_{11} & b_1A_{12} & \cdots & b_1\cdots b_{n-1}A_{1n}\\ & A_{22} & \ddots & \vdots\\ & & \ddots & b_{n-1}A_{n-1, \, n}\\ 0 & & & A_{nn}\end{array}\right]).$$

\vspace{.5cm}

{\em Proof}. \cite[Lemma 3.2]{11} says that if $|a|\le |b|$, then
$$W(\left[\begin{array}{cc} A & aB\\ 0  & C\end{array}\right])\subseteq W(\left[\begin{array}{cc} A & bB \\ 0 & C\end{array}\right]).$$
We apply this result $n-1$ times to obtain
\begin{eqnarray*}
&& W(\left[\begin{array}{c|cccc} A_{11} & a_1A_{12} & a_1a_2A_{13} & \cdots & a_1\cdots a_{n-1}A_{1n}\\ \hline & A_{22} & a_2A_{23} & \cdots & a_2\cdots a_{n-1}A_{2n}\\ & & \ddots & \ddots & \vdots\\ & & & \ddots & a_{n-1}A_{n-1, \, n}\\ & & & & A_{nn}\end{array}\right])\\
& \subseteq & W(\left[\begin{array}{c|cccc} A_{11} & b_1A_{12} & b_1a_2A_{13} & \cdots & b_1a_2\cdots a_{n-1}A_{1n}\\ \hline & A_{22} & a_2A_{23} & \cdots & a_2\cdots a_{n-1}A_{2n}\\ & & \ddots & \ddots & \vdots\\ & & & \ddots & a_{n-1}A_{n-1, \, n}\\ & & & & A_{nn}\end{array}\right])\\
& \subseteq & W(\left[\begin{array}{cc|cccc} A_{11} & b_1A_{12} & b_1b_2A_{13}& b_1b_2a_3A_{14} & \cdots & b_1b_2a_3\cdots a_{n-1}A_{1n}\\  & A_{22} & b_2A_{23} & b_2a_3A_{24} & \cdots & b_2a_3\cdots a_{n-1}A_{2n}\\ \hline & & A_{33} & a_3A_{34} & \cdots & a_3\cdots a_{n-1}A_{3n}\\ & & &\ddots & \ddots & \vdots\\ & & & & \ddots & a_{n-1}A_{n-1, \, n}\\ & & & & & A_{nn}\end{array}\right])\\
& \subseteq & \cdots\\
& \subseteq & W(\left[\begin{array}{cccc} A_{11} & b_1A_{12} & \cdots & b_1\cdots b_{n-1}A_{1n}\\ & A_{22} & \ddots & \vdots\\ & & \ddots & b_{n-1}A_{n-1, \, n}\\ 0 & & & A_{nn}\end{array}\right]). \hspace{60mm} \blacksquare
\end{eqnarray*}

\vspace{5mm}

{\em Proof of Proposition} 2.1. (a) If $a=e^{i\theta}b$ ($\theta\in\mathbb{R}$) and $U=\dia(1, e^{i\theta}, e^{2i\theta}, \ldots, e^{(n-1)i\theta})$, then $U$ is unitary and $UJ_n(a)=J_n(b)U$.

\vspace{3mm}

(b) If $X=\dia(1, a/b, (a/b)^2, \ldots, (a/b)^{n-1})$, then $X$ is invertible and $XJ_n(a)=J_n(b)X$.

\vspace{3mm}

(c) Since $J_n(a)$ is unitarily similar to $J_n(|a|)$ by (a), our assertion follows from the fact that the numerical range of a real matrix is always symmetric with respect to the $x$-axis.

\vspace{3mm}

(d) For $n\ge 2$ and $a\neq 0$, $J_n(a)$ is a nonzero nilpotent matrix with $J_n(a)^{n-1}\neq 0_n$. Hence $J_n(a)$ is irreducible (cf. proof of \cite[Theorem 3.1]{24}). The other two assertions follow from \cite[Corollary 1.2]{10}.

\vspace{3mm}

(e) This is obtained from Lemma 2.2 by letting $a_j=a$ and $b_j=b$ for all $j$, and $A_{ij}=0_1$ if $i=j$, and $I_1$ if $i<j$. \hfill $\blacksquare$

\vspace{5mm}

The next theorem characterizes those $J_n(a)$'s whose numerical ranges are circular discs. We will give two different proofs. For the first one, we need the Kippenhahn polynomial of a matrix. Recall that the \emph{Kippenhahn polynomial} of an $n$-by-$n$ matrix $A$ is the degree-$n$ real-coefficient homogeneous polynomial $p_A(x,y,z)$ given by $\det(x\re A+y\im A+zI_n)$. It relates to the numerical range of $A$ by the fact that $W(A)$ equals the convex hull of the real points of the dual curve of $p_A(x,y,z)=0$ (cf. \cite[Theorem 10]{17}).

\vspace{5mm}

{\bf Theorem 2.3.} \emph{The following statements are equivalent for} $J_n(a)$:

(a) $W(J_n(a))$ \emph{is a circular disc},

(b) \emph{the boundary of $W(J_n(a))$ contains an elliptic arc}, \emph{and}

(c) $n=2$ \emph{and} $a\neq 0$.

\vspace{5mm}

{\em Proof} 1. Obviously, (a) implies (b). To prove that (b) implies (c), assume that $n\ge 3$, $a\neq 0$, and $E$ is an elliptic disc such that $\partial W(J_n(a))$ contains an arc of $\partial E$. By Proposition 2.1 (a), we may further assume that $a>0$. Let $A=J_n(a)$ and let $B$ be a 2-by-2 matrix with $W(B)=E$. Via duality and B\'{e}zout's theorem \cite[Theorem 3.9]{18}, we infer that $p_B$ is a factor of $p_A$. In particular, $p_B(1, i, z)=\det(B+zI_2)$ divides $p_A(1, i,z)=\det(A+zI_n)$. Hence the two eigenvalues of $B$ are also eigenvalues of $A$. Thus they are both 0. Therefore, we may assume that $B={\scriptsize\left[\begin{array}{cc} 0 & 2b\\ 0 & 0\end{array}\right]}$ for some $b>0$. Simple computations then yield that $p_B(x,y,z)=z^2-b^2(x^2+y^2)$. On the other hand, we also have
\begin{eqnarray*}
&& p_A(1, y, 0)=\det(\re A+y\im A)\\[3mm]
&=& \det\Bigg(\frac{1}{2}(1-iy)\left[\begin{array}{cccc} 0 & a & \cdots & a^{n-1}\\ & 0 & \ddots & \vdots\\ & & \ddots & a\\ 0 & & & 0\end{array}\right]+\frac{1}{2}(1+iy)\left[\begin{array}{cccc} 0 &   &   & 0\\ a & 0 &   &  \\ \vdots & \ddots & \ddots &  \\ a^{n-1} & \cdots & a & 0\end{array}\right]\Bigg)\\[3mm]
&=& (\frac{a}{2})^n\det\left[\begin{array}{cccccc}
0 & 1-iy & a(1-iy) & \cdots & \cdots & a^{n-2}(1-iy)\\
1+iy & 0 & 1-iy & \ddots & & \vdots\\
a(1+iy) & 1+iy & 0 & \ddots & \ddots & \vdots\\
\vdots & \ddots & \ddots & \ddots & \ddots & a(1-iy)\\
\vdots & & \ddots & \ddots & \ddots & 1-iy\\
a^{n-2}(1+iy) & \cdots & \cdots & a(1+iy) & 1+iy & 0\end{array}\right]\\[3mm]
&=& (\frac{a}{2})^n\sum_{\scriptstyle 1\le j, k\le n-1\atop \scriptstyle j+k=n}\alpha_{jk}a^{m(j, \, k)}(1-iy)^j(1+iy)^k\\[3mm]
&=& (\frac{a}{2})^n\left[(-1)^{n-1}a^{n-2}(1+y^2)(1+iy)^{n-2}+\sum_{\scriptstyle 2\le j\le k\le n-2\atop \scriptstyle j+k=n}\alpha_{jk}a^{m(j, \, k)}(1+y^2)^j(1+iy)^{k-j}\right.\\[3mm]
&& \left.\mbox{\hspace{11mm}} +(-1)^{n-1}a^{n-2}(1-iy)^{n-2}(1+y^2)+\sum_{\scriptstyle 2\le k<j\le n-2\atop \scriptstyle j+k=n}\alpha_{jk}a^{m(j, \, k)}(1-iy)^{j-k}(1+y^2)^k\right],
\end{eqnarray*}
where $\alpha_{jk}$ can be $\pm 1$ or 0 and $m(j, k)$ is a nonnegative integer, both depending on the values of $j$ and $k$. Since $p_A=p_Bq$ for some degree-$(n-2)$ homogeneous polynomial $q$, we have
\begin{eqnarray*}
&& q(1, y, 0)=\frac{p_A(1,y,0)}{p_B(1,y,0)}=-\frac{1}{b^2}\frac{p_A(1,y,0)}{1+y^2}\\
&=& -\frac{1}{b^2}(\frac{a}{2})^n\left[(-1)^{n-1}a^{n-2}(1+iy)^{n-2}+\sum_{\scriptstyle 2\le j\le k\le n-2\atop \scriptstyle j+k=n}\alpha_{jk}a^{m(j, \, k)}(1+y^2)^{j-1}(1+iy)^{k-j}\right.\\
&& \left.\mbox{\hspace{18mm}} +(-1)^{n-1}a^{n-2}(1-iy)^{n-2}+\sum_{\scriptstyle 2\le k<j\le n-2\atop \scriptstyle j+k=n}\alpha_{jk}a^{m(j, \, k)}(1-iy)^{j-k}(1+y^2)^{k-1}\right].
\end{eqnarray*}
Plugging in $y=i$, we obtain
\begin{eqnarray}
&& q(1,i,0)=(-\frac{1}{b^2})(\frac{a}{2})^n(-1)^{n-1}a^{n-2}2^{n-2}\nonumber\\
&=& (-1)^n\frac{1}{4b^2}a^{2n-2}.\label{eq1}
\end{eqnarray}
Since $p_A(1,i,z)=z^n$, $p_B(1,i,z)=z^2$ and $p_A(1,i,z)=p_B(1,i,z)q(1,i,z)$, we obtain $q(1,i,z)=z^{n-2}$ and hence $q(1,i,0)=0$. We infer from (\ref{eq1}) that $a=0$, which is a contradiction. This proves (b) $\Rightarrow$ (c).

\vspace{3mm}

For (c) $\Rightarrow$ (a), note that
$$W(J_2(a))=W(\left[\begin{array}{cc} 0 & a\\ 0 & 0\end{array}\right])=\{z\in \mathbb{C}: |z|\le\frac{|a|}{2}\},$$
which yields (a). \hfill $\blacksquare$

\vspace{5mm}

We now proceed to the second proof of Theorem 2.3, which is based on the theories of $S_n$-, $S_n^{-1}$- and nonnegative matrices. Although the arguments go through more detours than the first proof, the auxiliary results obtained along the way are interesting on their own and should be useful for other occasions. We start with the $S_n$-matrices. An $n$-by-$n$ matrix $A$ is of \emph{class} $S_n$ if it is a \emph{contraction} ($\|A\|\equiv\max_{\|x\|=1}\|Ax\|\le 1$), its eigenvalues are all in $\mathbb{D}\equiv\{z\in \mathbb{C}: |z|<1\}$, and it satisfies $\rank(I_n-A^*A)=1$. Such matrices are the finite-dimensional version of the $S(\phi)$ operators ($\phi$ an inner function), which were first studied by Sarason \cite{21} in 1967 and featured prominently in the Sz.-Nagy--Foia\c{s} contraction theory \cite{22, 1}. On the other hand, $A$ is of \emph{class} $S_n^{-1}$ if all its eigenvalues have moduli greater than 1 and it satisfies $\rank (I_n-A^*A)=1$. Such matrices were first defined and studied by the first author \cite{4}. The KMS matrices $J_n(a)$ are related to these two classes of matrices via affine functions.

\vspace{5mm}

{\bf Lemma 2.4.} \emph{If} $0<|a|<1$ (\emph{respectively}, $|a|>1$), \emph{then $((1-|a|^2)/a)J_n(a)-\overline{a}I_n$ is of class $S_n$} (\emph{respectively}, \emph{of class} $S_n^{-1}$) \emph{with spectrum the singleton} $\{-\overline{a}\}$.

\vspace{5mm}

{\em Proof}. Let $A=((1-|a|^2)/a)J_n(a)-\overline{a}I_n$. A simple computation shows that
$$A=\left[\begin{array}{cccccc}
-\overline{a} & 1-|a|^2 & a(1-|a|^2) & \cdots & \cdots & a^{n-2}(1-|a|^2)\\
 & -\overline{a} & 1-|a|^2 & \ddots & & \vdots\\
 & & -\overline{a} & \ddots & \ddots & \vdots\\
 & & & \ddots & \ddots & a(1-|a|^2)\\
 & & & & \ddots & 1-|a|^2\\
0 & & & & & -\overline{a}
\end{array}\right].$$
Thus if $0<|a|<1$, then $A$ is of class $S_n$ by \cite[Theorem 1.2]{7}. On the other hand, if $|a|>1$, then $-A$ is of class $S_n^{-1}$ by \cite[Theorem 2.4]{4} and hence the same is true for $A$. \hfill $\blacksquare$

\vspace{5mm}

The following corollary is an easy consequence of Theorem 2.3 and Lemma 2.4.

\vspace{5mm}

{\bf Corollary 2.5.} \emph{Let $A$ be an $S_n$-matrix} (\emph{respectively}, \emph{$S_n^{-1}$-matrix}) \emph{with $\sigma(A)$ a singleton $\{\lambda\}$}. \emph{Then the following statements are equivalent}:

(a) \emph{$W(A)$ is a circular disc},

(b) \emph{$\partial W(A)$ contains an elliptic arc}, \emph{and}

(c) \emph{$n=2$ or $\lambda=0$} (\emph{respectively}, $n=2$).

\vspace{5mm}

{\em Proof}. We need only prove that (b) $\Rightarrow$ (c) and (c) $\Rightarrow$ (a). For the former, we assume that $\lambda\neq 0$. Then $(-\overline{\lambda}/(1-|\lambda|^2))(A-\lambda I_n)$ is unitarily similar to $J_n(-\overline{\lambda})$ by Lemma 2.4 and $\partial W(J_n(-\overline{\lambda}))$ contains an elliptic arc. Thus Theorem 2.3 implies that $n=2$. For (c) $\Rightarrow$ (a), if $n=2$, then $W(A)$ is an elliptic disc with foci both equal to $\lambda$, that is, it is a circular disc centered at $\lambda$. On the other hand, if $\lambda=0$, then $A$ is unitarily similar to the $n$-by-$n$ Jordan block
$$J_n=\left[\begin{array}{cccc} 0 & 1 &   & 0\\ & 0 & \ddots & \\ & & \ddots & 1\\ 0 & & & 0\end{array}\right]$$
(cf. \cite[Theorem 1.2]{7}). The numerical range of the latter equals $\{z\in \mathbb{C}: |z|\le\cos(\pi/(n+1))\}$ (cf. \cite[Proposition 1]{13}). This proves that (c) $\Rightarrow$ (a). \hfill $\blacksquare$

\vspace{5mm}

{\bf Lemma 2.6.} \emph{Let} $n\ge 3$, $a\neq 0$, \emph{and $\lambda$ be a point in $\partial W(J_n(a))$}. \emph{Then $\dim\bigvee\{x\in \mathbb{C}^n: \langle J_n(a)x, x\rangle=\lambda\|x\|^2\}=1$ if and only if either $|a|\neq 1$ or $|a|=1$ and $\re\lambda\neq -1/2$}.

\vspace{5mm}

{\em Proof}. If $|a|=1$, then $J_n(a)$ is unitarily similar to $J_n(1)$ by Proposition 2.1 (a). Properties of the numerical range of the latter were given in \cite[Lemma 1.3]{10}. In particular, if $\re\lambda=-1/2$, then $\partial W(J_n(1))$ has a vertical line segment passing through $\lambda$, which would imply that $\bigvee\{x\in \mathbb{C}^n: \langle J_n(1)x, x\rangle=\lambda\|x\|^2\}=\ker((\re J_n(1))-(-1/2)I_n)$ has dimension bigger than 1. This proves the necessity.

\vspace{3mm}

For the sufficiency, note that $A\equiv f(J_n(a))$, where $f(z)=((1-|a|^2)/a)z-\overline{a}$, is of class $S_n$ (respectively, of class $S_n^{-1}$) if $0<|a|<1$ (respectively, $|a|>1$) by Lemma 2.4, and $\eta\equiv f(\lambda)$ is in $\partial W(A)$. Hence \cite[Lemma 2.2]{5} (respectively, \cite[Theorem 2.5 (5)]{4}) yields that
$$\bigvee\{x\in \mathbb{C}^n: \langle Ax, x\rangle=\eta\|x\|^2\}=\bigvee\{x\in \mathbb{C}^n: \langle J_n(a)x, x\rangle=\lambda\|x\|^2\}$$
has dimension 1. On the other hand, if $|a|=1$ and $\re\lambda\neq -1/2$, then the dimension assertion follows from \cite[Lemma 1.3 (5)]{10}. \hfill $\blacksquare$

\vspace{5mm}

{\bf Lemma 2.7.} \emph{If $A$ is an $n$-by-$n$ matrix such that $\dim\bigvee\{x\in \mathbb{C}^n: \langle Ax, x\rangle=\lambda\|x\|^2\}=1$ for all $\lambda$ in $\partial W(A)$}, \emph{then there is an irreducible factor $q$ of $p_A$ such that $W(A)$ equals the convex hull of the real points of the dual curve of $q(x,y,z)=0$ and}, \emph{in particular}, \emph{$\partial W(A)$ is an irreducible algebraic curve}.

\vspace{5mm}

In general, if $A$ is any $n$-by-$n$ matrix with $p_A$ factored as the product $q_1\cdots q_m$ of irreducible factors, then, by \cite[Theorem 10]{17}, $W(A)$ is the convex hull of the union of the real points of the dual curves of $q_j(x,y,z)=0$, $1\le j\le m$. Hence, in particular, $\partial W(A)$ is the union of finitely many irreducible algebraic arcs.

\vspace{5mm}

{\em Proof of Lemma} 2.7. Let $p_A=q_1\cdots q_m$ with irreducible $q_j$'s, and let $\bigtriangleup_j$ be the convex hull of the real points of the dual of $q_j=0$, $1\le j\le m$. As noted above, we have $W(A)=(\cup_{j=1}^m\bigtriangleup_j)^{\wedge}$. Assume that $\lambda=a+bi$ ($a, b$ real) is a point in the intersection $\partial\bigtriangleup_j  \cap \, \partial\bigtriangleup_k  \cap \, \partial W(A)$ for some $j$ and $k$, $1\le j\neq k\le m$. If $L$ is the supporting line of $W(A)$ at $\lambda$, $R$ is the ray from the origin which is perpendicular to $L$, and $\theta$ is the angle from the positive $x$-axis to $R$, then $L$ is given by $x\cos\theta+y\sin\theta=a\cos\theta+b\sin\theta=\re(e^{-i\theta}\lambda)\equiv d$ and is also a supporting line of both $\bigtriangleup_j$ and $\bigtriangleup_k$ at $\lambda$. By duality, we have $q_j(\cos\theta, \sin\theta, -d)=q_k(\cos\theta, \sin\theta, -d)=0$. Hence $-d$ is a zero of the polynomial $p_A(\cos\theta, \sin\theta,z)$ in $z$ with multiplicity at least 2. Since $p_A(\cos\theta, \sin\theta,z)=\det(\re(e^{-i\theta}A)+zI_n)$, this means that $d$ is an eigenvalue of $\re(e^{-i\theta}A)$ of multiplicity at least 2. Note that its corresponding eigenspace $\{x\in \mathbb{C}^n: \re(e^{-i\theta}A)x=dx\}$ is the same as $\bigvee\{x\in \mathbb{C}^n: \langle Ax, x\rangle=\lambda\|x\|^2\}$ with dimension 1 by our assumption. Thus the contradiction yields that $\partial\bigtriangleup_j \cap \, \partial\bigtriangleup_k \cap \, \partial W(A)=\emptyset$ for any $j\neq k$. Since $\partial W(A)$ contains no line segment by our assumption, we conclude that $W(A)=\bigtriangleup_j$ for some $j$ and, in particular, $\partial W(A)$ is an irreducible algebraic curve. \hspace{-3mm} \hfill $\blacksquare$

\vspace{5mm}

{\bf Lemma 2.8.} \emph{Let $A$ be an $n$-by-$n$ matrix such that $W(A)$ is the convex hull of an irreducible algebraic curve}. \emph{If $\partial W(A)$ contains an arc of the ellipse $C$}, \emph{then} $W(A)=C^{\wedge}$.

\vspace{5mm}

{\em Proof}. Let $\partial W(A)$ (respectively, $C$) be given by the irreducible algebraic (respectively, quadratic) curve $q_1=0$ (respectively, $q_2=0$) together with possibly some line segments. Our assumption on $\partial W(A)$ and $C$ implies, by duality and B\'{e}zout's theorem \cite[Theorem 3.9]{18}, that $q_1=q_2$. Hence $W(A)=C^{\wedge}$ follows. \hfill $\blacksquare$

\vspace{5mm}

The next corollary shows that the preceding lemmas are applicable to $S_n$- and $S_n^{-1}$-matrices. It answers a question asked by K.-Z. Wang.

\vspace{5mm}

{\bf Corollary 2.9.} \emph{Let $A$ be an $S_n$-matrix} (\emph{respectively}, \emph{$S_n^{-1}$-matrix}). \emph{If $\partial W(A)$ contains an arc of the ellipse $C$}, \emph{then} $W(A)=C^{\wedge}$.

\vspace{5mm}

{\em Proof}. Since $A$ satisfies $\dim\bigvee\{x\in \mathbb{C}^n: \langle Ax, x\rangle=\lambda\|x\|^2\}=1$ for all $\lambda$ in $\partial W(A)$ by \cite[Lemma 2.2]{5} (respectively, \cite[Theorem 2.5 (5)]{4}), our assertion follows from Lemmas 2.7 and 2.8. \hfill $\blacksquare$

\vspace{5mm}

Note that if $A$ is of class $S_n^{-1}$, then the preceding corollary also follows directly from Lemma 2.8 since, in this case, $p_A$ is irreducible (cf. \cite[Lemma 2.9 (2)]{4}).

\vspace{5mm}

We next move to properties of nonnegative matrices. recall that a matrix $A=[a_{ij}]_{i, j=1}^n$ is \emph{nonnegative}, denoted by $A\succcurlyeq 0$, if $a_{ij}\ge 0$ for all $i$ and $j$. It is \emph{permutationally irreducible} if there is no \emph{permutation matrix} $P$ (a matrix with every row and column having exactly one entry 1 and all others 0) such that $P^TAP$ is of the form ${\scriptsize\left[\begin{array}{cc} B & C\\ 0 & D\end{array}\right]}$, where $B$ and $D$ are square matrices. Properties of the numerical ranges of nonnegative matrices were given in \cite{19}.

\vspace{5mm}

{\bf Proposition 2.10.} \emph{If} $n\ge 3$, \emph{then the boundary of $W(J_n(a))$ intersects the circle $|z|=w(J_n(a))$ at exactly one point}, \emph{namely}, \emph{at} $w(J_n(a))$.

\vspace{5mm}

{\em Proof}. By Proposition 2.1 (a), we may assume that $a>0$. Then $A\equiv J_n(a)\succcurlyeq 0$ and $\re A$ is permutationally irreducible. Hence $w(A)$ is in $\partial W(A)\cap\{z\in\mathbb{C}:|z|=w(A)\}$ (cf. \cite[p. 5]{19}). We now show that it is the only point in this intersection. Indeed, assume that $e^{i\theta}w(A)$ is in $\partial W(A)$ for some real $\theta$. \cite[Proposition 3.7]{19} says that $A$ is unitarily similar to $e^{i\theta}A$. Hence
\begin{equation}\tr(A^2A^*)=\tr((e^{i\theta}A)^2(e^{i\theta}A)^*)=e^{i\theta}\tr(A^2A^*).\label{eq2}\end{equation}
Since
\begin{eqnarray*}
&& A^2A^*=\left[
            \begin{array}{cccccc}
              0 & 0 & a^2 & 2a^3 & \cdots & (n-2)a^{n-1} \\
                & 0 & 0 & \ddots & \ddots & \vdots \\
             &   & \ddots  & \ddots & \ddots & 2a^3 \\
              &   &   & \ddots & \ddots & a^2 \\
             &   &   &   & \ddots & 0 \\
             &    &   &   &  & 0
            \end{array}
          \right]
          \left[
            \begin{array}{cccc}
              0&   &   &   \\
              a & 0 &   &   \\
              \vdots & \ddots & \ddots&   \\
             a^{n-1} & \cdots & a & 0 \\
            \end{array}
          \right]\\
&=& \left[
      \begin{array}{cc}
        \sum_{k=1}^{n-2}ka^{2(k+1)} & \ast \\
        \ast & \ddots \\
      \end{array}
    \right]\succcurlyeq 0,
\end{eqnarray*}
we have $\tr(A^2A^*)\ge\sum_{k=1}^{n-2}ka^{2(k+1)}>0$, which, together with (\ref{eq2}), yields that $e^{i\theta}=1$. This proves that $\partial W(A)\cap\{z\in\mathbb{C}: |z|=w(A)\}=\{w(A)\}$. \hfill $\blacksquare$

\vspace{5mm}

A combination of the previous lemmas with Proposition 2.10 yields the second proof of Theorem 2.3.

\vspace{5mm}

{\em Proof $2$ of Theorem $2.3$}. We need only prove (b) $\Rightarrow$ (c). Assume that (b) holds. Let $A=J_n(a)$ for $n\ge 3$ and $a>0$. If $a\neq 1$, then $\dim\bigvee\{x\in \mathbb{C}^n: \langle Ax, x\rangle=\lambda\|x\|^2\}=1$ for all $\lambda$ in $\partial W(A)$ by Lemma 2.6. We infer from Lemmas 2.7 and 2.8 that $W(A)$ equals an elliptic disc, say, $E$. Since the foci of $\partial E$ are the eigenvalues 0 of $A$ (cf. \cite[Theorem]{8}), $W(A)$ ($=E$) is a circular disc centered at the origin. This contradicts Proposition 2.10. On the other hand, if $a=1$, then $p_A$ is irreducible by \cite[Lemma 1.3 (8)]{10}. Again, we can infer from Lemma 2.8 that $W(A)=E$, which contradicts the fact that $\partial W(A)$ contains a line segment (cf. \cite[Lemma 1.3 (4)]{10}). This proves (c). \hspace{-2mm}$\blacksquare$

\vspace{5mm}

Before we move on, two remarks are in order. Firstly, in the preceding proof, after we have shown that $W(A)$ is a circular disc centered at the origin, we may resort to \cite[Theorem 3]{23}, instead of Proposition 2.10, to reach the contradiction. This is because one pair of the equivalent conditions there says that, for an $n$-by-$n$ nonnegative matrix $A$ with $\re A$ permutationally irreducible, $W(A)$ is a circular disc centered at the origin if and only if $A$ is permutationally similar to a matrix of the form
$$\left[
    \begin{array}{cccc}
      0 & A_1 &   & 0 \\
        & 0 & \ddots &   \\
        &   & \ddots & A_{m-1} \\
      0  &   &   & 0 \\
    \end{array}
  \right]$$
for some $m\ge 2$, and $J_n(a)$ obviously does not satisfy this condition for $n\ge 3$ and $a>0$. Secondly, \cite[Theorem 4.5]{20} on the noncircularity of $W(A)$ for a permutationally irreducible nonnegative matrix $A$ is not applicable in Proof 2 of Theorem 2.3 since $J_n(a)$ itself is not permutationally irreducible.

\vspace{5mm}

Another benefit of Lemma 2.4 is that we can use known properties of the numerical ranges of $S_n$- and $S_n^{-1}$-matrices to deduce the corresponding ones for $J_n(a)$. The following proposition is one such example.

\vspace{5mm}

{\bf Proposition 2.11.} \emph{The boundary of $W(J_n(a))$ contains a line segment if and only if $n\ge 3$ and $|a|=1$}.

\vspace{5mm}

{\em Proof}. If $\partial W(J_n(a))$ contains a line segment, then $n\ge 3$ and $a\neq 0$ trivially. Assume that $|a|\neq 1$. Then Lemma 2.4 says that $A\equiv f(J_n(a))$, where $f(z)=((1-|a|^2)/a)z-\overline{a}$, is in class $S_n$ or $S_n^{-1}$ depending on whether $0<|a|<1$ or $|a|>1$. Since the boundary of $W(A)$ contains no line segment by \cite[Lemma 2.2]{5} or \cite[Theorem 2.5 (4)]{4}, the same is true for the boundary of $W(J_n(a))$. This contradicts our assumption. Thus we must have $|a|=1$.

\vspace{3mm}

Conversely, if $n\ge 3$ and $|a|=1$, then $J_n(a)$ is unitarily similar to $J_n(1)$ by Proposition 2.1 (a). It is known from \cite[Lemma 1.3 (4)]{10} that $\partial W(J_n(1))$ contains a line segment. Thus the same is true for $W(J_n(a))$. \hfill $\blacksquare$

\vspace{5mm}

We conclude this section with some partial results on the irreducibility of the Kippenhahn polynomial of $J_n(a)$. From the results we obtained so far, it seems likely that $p_{J_n(a)}$ is irreducible for any $a\neq 0$. If this is indeed the case, then in Proof 2 of Theorem 2.3 we can bypass Lemmas 2.6 and 2.7 to obtain the noncircularity of $W(J_n(a))$ directly. However, at this point we are not able to confirm this. The following proposition contains the positive special cases which we know of.

\vspace{5mm}

{\bf Proposition 2.12.} \emph{If either} (a) $n\le 4$ \emph{and} $a\neq 0$, \emph{or} (b) $|a|>\cos(\pi/(n+1))$, \emph{then the Kippenhahn polynomial $p_{J_n(a)}$ is irreducible}.

\vspace{5mm}

{\em Proof}. Assume that $a>0$.

\vspace{3mm}

(a) If $p_{J_2(a)}$ is reducible, say, $p_{J_2(a)}(x,y,z)=(a_1x+b_1y+z)(a_2x+b_2y+z)$ for some real $a_j$ and $b_j$, $j=1, 2$, then $W(J_2(a))$ is the line segment connecting $a_1+b_1i$ and $a_2+b_2i$. This implies that $J_2(a)$ is a normal matrix, which is certainly not the case.

\vspace{3mm}

Next for $n=3$. We will make use of the classification of the numerical ranges of 3-by-3 matrices via their Kippenhahn polynomials from \cite[Theorem 2.6]{17} or \cite{16}. If $p_{J_3(a)}$ is the product of three linear factors, then, as for $n=2$, this corresponds to $J_3(a)$ being normal, again a contradiction. Next we assume that $p_{J_3(a)}=q_1q_2$, where $q_1$ is quadratic irreducible and $q_2$ is linear. Then $W(J_3(a))$ is an elliptic disc with possibly a cone attached to it. Either way, $\partial W(J_3(a))$ contains an elliptic arc. The equivalence of (b) and (c) of Theorem 2.3 yields that $n=2$, which is a contradiction. This proves the irreducibility of $p_{J_3(a)}$.

\vspace{3mm}

We now consider the case of $n=4$. Let $A=J_4(a)$ and assume that $p_A$ is reducible, say, $p_A=q_1q_2$. As for $n=3$, if $q_1$ and $q_2$ are both reducible, then $A$ is normal, a contradiction. On the other hand, if at least one of the $q_j$'s is quadratic irreducible, then $\partial W(A)$ contains an elliptic arc and then Theorem 2.3 leads to a contradiction as before. We are thus left with the case of a cubic irreducible $q_1$ and linear $q_2$. Let $q_2(x,y,z)=cx+dy+z$, where $c$ and $d$ are real. Since
$$p_A(1,i,z)=q_1(1,i,z)q_2(1,i,z)=q_1(1,i,z)(c+di+z)$$
and
$$p_A(1,i,z)=\det(A+zI_4)=z^4,$$
we infer that $c=d=0$ and $q_2(x,y,z)=z$. Hence $p_A(1,y,z)=q_1(1,y,z)z$. This shows that the term in the expansion of $p_A(1,y,z)$ which contains no $z$ must be 0. But
\begin{eqnarray*}
&& p_A(1,y,z)=\det(\re A+y\im A+zI_4)\\
&=& \left[
      \begin{array}{cccc}
        z & a(1-iy)/2 & a^2(1-iy)/2  & a^3(1-iy)/2  \\
        a(1+iy)/2  & z & a(1-iy)/2  & a^2(1-iy)/2  \\
         a^2(1+iy)/2 &  a(1+iy)/2 & z & a(1-iy)/2  \\
         a^3(1+iy)/2 &  a^2(1+iy)/2 &  a(1+iy)/2 & z \\
      \end{array}
    \right].
\end{eqnarray*}
The term referred to above can be calculated to be $a^4(1+y^2)(1-4a^2+y^2)/16$, which is nonzero. We conclude that $p_A$ is irreducible.

\vspace{3mm}

(b) Note that if $a>1$, then $A\equiv ((1-a^2)/a)J_n(a)-aI_n$ is of class $S_n^{-1}$ by Lemma 2.4. Since $p_A$ is irreducible by \cite[Lemma 2.9 (2)]{4}, the same is true for $p_{J_n(a)}$. For $a=1$, the irreducibility of $p_{J_n(1)}$ is proven in \cite[Lemma 1.3 (8)]{10}. Thus we need only consider $\cos(\pi/(n+1))<a<1$. For this case, let $\lambda_1\ge\lambda_2\ge\cdots\ge\lambda_n$ be the eigenvalues of $\re J_n(a)$. As observed in \cite[pp. 69--70]{12}, the $\lambda_k$'s are distinct and $2\lambda_k+1$, $1\le k\le n$, are eigenvalues of the matrix $[a^{|i-j|}]_{i, j=1}^n$, which satisfy $2\lambda_k+1=(1-a^2)/|1-ae^{it_k}|^2$, where $t_1<t_2<\cdots<t_n$ are the roots of $\sin((n+1)t)-2a\sin(nt)+a^2\sin((n-1)t)=0$. Moreover, it is known that $t_2>\pi/(n+1)$ (cf. \cite[p. 70]{12}). Thus we have
$$2\lambda_2+1=\frac{1-a^2}{|1-ae^{it_2}|^2}<\frac{1-a^2}{|1-ae^{i(\pi/(n+1))}|^2}=\frac{1-a^2}{1-2a\cos(\pi/(n+1))+a^2}<1,$$
which implies that $\lambda_2<0$. Now assume that $p_{J_n(a)}=q_1q_2$ with $q_1(1, 0, -\lambda_1)=0$. Let $r$ be the largest root of $q_2(1, 0, -z)=0$. Then $0=q_1(1, 0, -r)q_2(1, 0 -r)=p_{J_n(a)}(1, 0, -r)=\det(\re J_n(a) -rI_n)$. Since the $\lambda_k$'s are distinct, this guarantees that $r=\lambda_{k_0}$ for some $k_0\ge 2$. Let $\bigtriangleup$ be the convex hull of the real points of the dual curve of $q_2=0$. We infer by duality that the vertical line $x=\lambda_{k_0}$ is a supporting line of $\bigtriangleup$ or, in other words, $\re\bigtriangleup\le\lambda_{k_0}\le\lambda_2<0$. This means that $\bigtriangleup$ is contained in the open left half-plane. On the other hand, the real foci of $\partial\bigtriangleup$, that is, the points $c+di$ ($c, d$ real) satisfying $q_2(1, \pm i, -(c\pm di))=0$ are all in $\bigtriangleup$ (cf. \cite[Lemma 2.8]{4}). Such points also satisfy $p_{J_n(a)}(1, \pm i, -(c\pm di))=0$ and hence must all be 0. These show that 0 is in $\bigtriangleup$, which contradicts what was proven before. Hence $p_{J_n(a)}$ must be irreducible. \hfill $\blacksquare$

\vspace{1cm}

\noindent
{\bf\large 3. Compressions}

\vspace{3mm}

In this section, we are concerned with the containment relations between the numerical ranges of a KMS matrix and its compressions. We start with the following proposition.

\vspace{5mm}

{\bf Proposition 3.1.} \emph{If} $1\le m<n$, $a\neq 0$, \emph{and $A$ is an $m$-by-$m$ compression of $J_n(a)$}, \emph{then} $W(A)\subsetneqq W(J_n(a))$.

\vspace{5mm}

{\em Proof}. Let $B=f(J_n(a))$, where $f(z)=((1-|a|^2)/a)z-\overline{a}$. If $|a|<1$, then $B$ is of class $S_n$ by Lemma 2.4 and $A'\equiv f(A)$ is a compression of $B$. Since $W(A')\subsetneqq W(B)$ by \cite[Theorem 3.3 (a)]{11}, we have $W(A)\subsetneqq W(J_n(a))$. On the other hand, if $|a|\ge 1$, then, by Proposition 2.12 (b), $p_{J_n(a)}$ is irreducible. Note that if $W(A)$ and $W(J_n(a))$ are equal, then $p_A$ and $p_{J_n(a)}$ have a common irreducible factor (cf. \cite[Proposition 2.3]{9}). The irreducibility of $p_{J_n(a)}$ yields that $p_{J_n(a)}$ divides $p_A$. Thus we have $n\le m$, which contradicts our assumption. Hence, in this case, we also have $W(A)\subsetneqq W(J_n(a))$. \hfill $\blacksquare$

\vspace{5mm}

A property closely related to the nonequality of the numerical ranges of a matrix and its compressions is given in the following.

\vspace{5mm}

{\bf Lemma 3.2.} \emph{If $A$ is an $n$-by-$n$ matrix with the property that $W(B)\subsetneqq W(A)$ for all $m$-by-$m$} ($1\le m<n$) \emph{compressions $B$ of} $A$, \emph{then} $\bigvee\{x\in\mathbb{C}^n: \|x\|=1, \langle Ax, x\rangle\in\partial W(A)\}=\mathbb{C}^n$. \emph{If}, \emph{moreover}, $A$ \emph{satisfies $\dim\bigvee\{x\in\mathbb{C}^n: \langle Ax, x\rangle=\lambda\|x\|^2\}=1$ for all $\lambda$ in $\partial W(A)$}, \emph{then the converse also holds}.

\vspace{5mm}

Note that, in the preceding lemma, the converse fails if there is no extra condition on $A$, as witness $A={\scriptsize\left[\begin{array}{cc} 0 & 0\\ 0 & 0\end{array}\right]}$ and $B=[0]$.

\vspace{5mm}

{\em Proof of Lemma $3.2$}. Let $K=\bigvee\{x\in\mathbb{C}^n: \|x\|=1, \langle Ax, x\rangle\in\partial W(A)\}$, $m=\dim K$, and $B$ be the compression of $A$ to $K$. Then, for any point $\lambda$ in $\partial W(A)$, there is a unit vector $x$ in $K$ such that $\langle Ax, x\rangle=\lambda$. Since $\langle Bx, x\rangle=\langle Ax, x\rangle=\lambda$, this shows that $\partial W(A)\subseteq W(B)$. We infer from the convexity of $W(A)$ and $W(B)$ that $W(A)\subseteq W(B)$. Since $W(B)\subseteq W(A)$ always holds, we obtain $W(A)=W(B)$. Hence our assumption on $A$ yields that $m=n$ or $K=\mathbb{C}^n$ as required.

\vspace{3mm}

Under the extra condition on the boundary points of $W(A)$, the converse was proven in \cite[Lemma 3.1 (a)]{11}. \hfill $\blacksquare$

\vspace{5mm}

The next corollary is a consequence of the preceding two results.

\vspace{5mm}

{\bf Corollary 3.3.} \emph{For any $n\ge 1$ and $a$ in $\mathbb{C}$}, \emph{the equality $\bigvee\{x\in\mathbb{C}^n: \|x\|=1, \langle J_n(a)x, x\rangle\in\partial W(J_n(a))\}=\mathbb{C}^n$ holds}.

\vspace{5mm}

Our next goal is to show that, for $n\ge 2$ and $|a|\neq 0, 1$, the numerical range of the restriction of $J_n(a)$ to one of its invariant subspaces is even contained in the interior of $W(J_n(a))$. This will be proven as a consequence of the following lemma.

\vspace{5mm}

{\bf Lemma 3.4.} \emph{If $a\neq 0$ and $x$ is a unit vector in $\mathbb{C}^n$ such that $\langle J_n(a)x, x\rangle$ is in the boundary of $W(J_n(a))$ but not in its boundary line segment}, \emph{then $x$ is a cyclic vector for $J_n(a)$}.

\vspace{5mm}

Recall that a vector $x$ in $\mathbb{C}^n$ is \emph{cyclic} for the $n$-by-$n$ matrix $A$ if $\bigvee\{x, Ax, \ldots, A^{n-1}x\}=\mathbb{C}^n$.

\vspace{5mm}

{\em Proof of Lemma $3.4$}. Assume that $n\ge 2$ and $a>0$, and let $\lambda=\langle J_n(a)x, x\rangle$. If $0<a<1$ (respectively, $a>1$), then $A\equiv f(J_n(a))$, where $f(z)=((1-a^2)/a)z-a$, is of class $S_n$ (respectively, of class $S_n^{-1}$) by Lemma 2.4, and $\langle Ax, x\rangle=f(\lambda)$ is in $\partial W(A)$. Hence $x$ is a cyclic vector for $A$ by \cite[Lemma 3.2]{6} (respectively, \cite[Theorem 2.5 (3)]{4}). Thus $x$ is also cyclic for $J_n(a)$.

\vspace{3mm}

Now assume that $a=1$ and let $x=[x_1 \ \ldots \ x_n]^T$. Our assumption on $\lambda$ implies that either $n=2$ or $n\ge 3$ and $\re\lambda\neq -1/2$ (cf. \cite[Lemma 1.3 (4)]{10}). In either case, we claim that $x_n\neq 0$. Indeed, if $n=2$ and $x_2=0$, then
$$\lambda=\langle J_2(1)x, x\rangle=\langle 0, x\rangle=0,$$
contradicting the fact that $W(J_2(1))=\{z\in \mathbb{C}: |z|\le 1/2\}$. On the other hand, if $n\ge 3$ and $\re\lambda\neq -1/2$, then $x_n\neq 0$ follows from the proof of \cite[Lemma 1.3 (4)]{10}. Since $J_n(1)^jx=[\ \ldots \ x_n \ \underbrace{0 \ \ldots \ 0}_j]^T$ for $1\le j\le n-1$, $x_n\neq 0$ implies that $\mathbb{C}^n$ is spanned by $x, J_n(1)x, \ldots, J_n(1)^{n-1}x$, that is, $x$ is cyclic for $J_n(1)$. \hfill $\blacksquare$

\vspace{5mm}

{\bf Proposition 3.5.} \emph{If} $n\ge 2$, $|a|\neq 0, 1$, \emph{and $J_n(a)$ is unitarily similar to} ${\scriptsize\left[\begin{array}{cc} A & \ast\\ 0 & \ast\end{array}\right]}$, \emph{where $A$ is of size $m$} ($1\le m<n$), \emph{then $W(A)$ is contained in the interior of $W(J_n(a))$}.

\vspace{5mm}

{\em Proof}. If $\lambda$ is any point in $\partial W(A)\cap\partial W(J_n(a))$, then $\lambda=\langle Ax, x\rangle$ for some unit vector $x$ in $\mathbb{C}^m$. Let $A'={\scriptsize\left[\begin{array}{cc} A & \ast\\ 0 & \ast\end{array}\right]}$ and $U$ be an $n$-by-$n$ unitary matrix such that $UJ_n(a)=A'U$. Since
$$\langle J_n(a)U^*(x\oplus 0), U^*(x\oplus 0)\rangle=\langle\left[\begin{array}{cc} A & \ast\\ 0 & \ast\end{array}\right]\left[\begin{array}{c}  x\\ 0 \end{array}\right], \left[\begin{array}{c}  x\\ 0 \end{array}\right]\rangle=\langle Ax, x\rangle=\lambda,$$
the unit vector $U^*(x\oplus 0)$ is cyclic for $J_n(a)$ (cf. Proposition 2.11 and Lemma 3.4). This implies that $x\oplus 0$ is cyclic for $A'$. But $\bigvee\{x\oplus 0, A'(x\oplus 0), \ldots, {A'}^{n-1}(x\oplus 0)\}\subseteq\mathbb{C}^m\oplus 0\neq \mathbb{C}^n$. The contradiction yields that $\partial W(A)\cap\partial W(J_n(a))=\emptyset$. Since $W(A)\subseteq W(J_n(a))$, the disjointness of $\partial W(A)$ and $\partial W(J_n(a))$ is equivalent to $W(A)$ being contained in the interior of $W(J_n(a))$. \hfill $\blacksquare$

\vspace{5mm}

{\bf Corollary 3.6.} \emph{If $1\le m<n$ and $|a|\neq 0, 1$}, \emph{then $W(J_m(a))$ is contained in the interior of $W(J_n(a))$}.

\vspace{5mm}

{\em Proof}. Since
$$J_n(a)=\left[\begin{array}{c|c} J_m(a) & \ast\\ \hline  0 & \begin{array}{cccc} 0 & a & \cdots & a^{n-m-1}\\ & 0 & \ddots & \vdots\\ & & \ddots & a\\ 0 & & & 0\end{array}\end{array}\right],$$
our assertion follows from Proposition 3.5. \hfill $\blacksquare$

\vspace{5mm}

Using the preceding corollary, we can refine the assertion in Proposition 2.1 (e).

\vspace{5mm}

{\bf Proposition 3.7.} \emph{If $n\ge 2$ and $|a|<|b|$}, \emph{then $W(J_n(a))$ is contained in the interior of $W(J_n(b))$}.

\vspace{5mm}

{\em Proof}. By Proposition 2.1 (a) and (d), we may assume that $0<a<b$. Let
$$A=\left[\begin{array}{cc} 0 & by^T\\ 0 & J_{n-1}(a)\end{array}\right],$$
where $y^T=[1 \ a \ \ldots \ a^{n-2}]$. Following the proof of Lemma 2.2, we have $W(J_n(a))\subseteq W(A)\subseteq W(J_n(b))$. Thus it suffices to show that $W(J_n(a))$ is contained in the interior of $W(A)$ or, equivalently, $\max\sigma(\re (e^{i\theta}J_n(a)))<\max\sigma(\re (e^{i\theta}A))$ for all real $\theta$. For this, we fix a $\theta$ in $\mathbb{R}$, and let $r=\max\sigma(\re (e^{i\theta}J_n(a)))$ and $x=x_1\oplus x_2$ in $\mathbb{C}^n=\mathbb{C}\oplus \mathbb{C}^{n-1}$ be a unit vector such that $(\re (e^{i\theta}J_n(a)))x=rx$. Then
\begin{eqnarray}
&& r=\langle (\re (e^{i\theta}J_n(a)))x, x\rangle \nonumber\\[3mm]
&=&\langle \left[\begin{array}{cc} 0 & (e^{i\theta}/2)ay^T\\ (e^{-i\theta}/2)ay & \re (e^{i\theta}J_{n-1}(a))\end{array}\right]\left[\begin{array}{c} x_1\\ x_2\end{array}\right], \left[\begin{array}{c} x_1\\ x_2\end{array}\right]\rangle \nonumber\\[3mm]
&=& \frac{e^{i\theta}}{2}ay^Tx_2\overline{x}_1+\frac{e^{-i\theta}}{2}ax_1\langle y, x_2\rangle+\langle (\re (e^{i\theta}J_{n-1}(a)))x_2, x_2\rangle \nonumber\\[3mm]
&=& a\re (e^{i\theta}\overline{x}_1\langle x_2, y\rangle)+\langle (\re (e^{i\theta}J_{n-1}(a)))x_2, x_2\rangle \label{e3}\\[3mm]
&\le & a|e^{i\theta}\overline{x}_1\langle x_2, y\rangle|+\langle (\re (e^{i\theta}J_{n-1}(a)))x_2, x_2\rangle \nonumber\\[3mm]
&=& \langle \left[\begin{array}{cc} 0 & (e^{i\theta}/2)ay^T\\ (e^{-i\theta}/2)ay & \re (e^{i\theta}J_{n-1}(a))\end{array}\right]\left[\begin{array}{c} x_1'\\ x_2\end{array}\right], \left[\begin{array}{c} x_1'\\ x_2\end{array}\right]\rangle \nonumber\\[3mm]
&\le & \max W(\re (e^{i\theta}J_{n}(a)))=r, \nonumber
\end{eqnarray}
where $x_1'=x_1e^{i\alpha}$ ($\alpha\in \mathbb{R}$) is such that $\re(ze^{-i\alpha})=|z|$ for $z=e^{i\theta}\overline{x}_1\langle x_2, y\rangle$. We infer from the above that $e^{i\theta}\overline{x}_1\langle x_2, y\rangle=|e^{i\theta}\overline{x}_1\langle x_2, y\rangle|\ge 0$. Furthermore, we claim that $e^{i\theta}\overline{x}_1\langle x_2, y\rangle>0$. Indeed, if $e^{i\theta}\overline{x}_1\langle x_2, y\rangle=0$, then either $x_1=0$ or $\langle x_2, y\rangle=0$ and we have $r=\langle (\re (e^{i\theta}J_{n-1}(a)))x_2, x_2\rangle$ by (\ref{e3}). If $x_1=0$, then $x_2$ is a unit vector in $\mathbb{C}^{n-1}$ and $r=\langle (\re (e^{i\theta}J_{n}(a)))x, x\rangle=\langle (\re (e^{i\theta}J_{n-1}(a)))x_2, x_2\rangle$ yields that $\partial W(J_n(a))\cap \partial W(J_{n-1}(a))\neq\emptyset$, which contradicts Corollary 3.6. On the other hand, if $\langle x_2, y\rangle=0$, then
\begin{eqnarray*}
&& r=\|x_2\|^2\langle (\re (e^{i\theta}J_{n-1}(a)))\frac{x_2}{\|x_2\|}, \frac{x_2}{\|x_2\|}\rangle\\[3mm]
&\le & \max W(\re (e^{i\theta}J_{n-1}(a))) < \max W(\re (e^{i\theta}J_{n}(a)))=r,
\end{eqnarray*}
where the last (strict) inequality is by Corollary 3.6, again a contradiction. Hence we derive from (\ref{e3}) that
\begin{eqnarray*}
&& r=ae^{i\theta}\overline{x}_1\langle x_2, y\rangle+\langle (\re (e^{i\theta}J_{n-1}(a)))x_2, x_2\rangle\\
&<& be^{i\theta}\overline{x}_1\langle x_2, y\rangle+\langle (\re (e^{i\theta}J_{n-1}(a)))x_2, x_2\rangle\\[3mm]
&=& \langle \left[\begin{array}{cc} 0 & (e^{i\theta}/2)by^T\\ (e^{-i\theta}/2)by & \re (e^{i\theta}J_{n-1}(a))\end{array}\right]\left[\begin{array}{c} x_1\\ x_2\end{array}\right], \left[\begin{array}{c} x_1\\ x_2\end{array}\right]\rangle \\[3mm]
&=& \langle (\re(e^{i\theta}A))x, x\rangle \le \max W(\re(e^{i\theta}A))\\
&=& \max\sigma(\re(e^{i\theta}A))
\end{eqnarray*}
as asserted. \hfill $\blacksquare$

\vspace{3mm}

Finally, we consider the relations between the numerical range of $J_n(a)$ and those of its principal submatrices. For any $n$-by-$n$ ($n\ge 2$) matrix $A$, we use $A[j]$, $1\le j\le n$, to denote its $j$th principal submatrix, that is, the $(n-1)$-by-$(n-1)$ matrix obtained by deleting the $j$th row and $j$th column of $A$.

\vspace{3mm}

{\bf Theorem 3.8.} \emph{For any $n\ge 2$ and $|a|\neq 0, 1$}, \emph{let} $b=\min\sigma(\re J_n(a))$. \emph{Then for any} $j$, $1\le j\le n$, \emph{we have}
$$\partial W(J_n(a))\cap\partial W(J_n(a)[j])=\left\{\begin{array}{ll} \{b\} & \mbox{\em if } \, n \ \mbox{\em is odd}, j=(n+1)/2, \mbox{\em and } |a|>1,\\ \emptyset & \mbox{\em otherwise}.\end{array}\right.$$
\emph{Moreover}, \emph{when $b$ is in this intersection}, \emph{there is a unit vector $x=[x_1 \ \ldots \ x_{j-1} \ 0 \ -x_{j-1} \ \ldots \ -x_1]^T$ in $\mathbb{C}^n$ such that $(\re J_n(a))x=bx$ and $\langle J_n(a)x, x\rangle=b$}.

\vspace{3mm}

Note that if $|a|=1$, then (a) it is easily seen that $\partial W(J_2(a))\cap\partial W(J_2(a)[j])=\emptyset$ for $j=1, 2$, (b) $\partial W(J_3(a))\cap\partial W(J_3(a)[j])=\{-1/2\}$ for $j=1, 2, 3$ by \cite[Lemma 1.3 (7)]{10}, and (c) for $n\ge 4$, $\partial W(J_n(a))\cap\partial W(J_n(a)[j])$ equals the line segment of $\partial W(J_{n-1}(1))$ on the vertical line $x=-1/2$ since $J_n(a)[j]$, $1\le j\le n$, is unitarily similar to $J_{n-1}(1)$ and $\partial W(J_m(1))$ contains a line segment on $x=-1/2$ for all $m\ge 3$ (cf. \cite[Lemma 1.3 (4)]{10}).

\vspace{3mm}

The proof of Theorem 3.8 will be done in a series of lemmas.

\vspace{3mm}

{\bf Lemma 3.9.} \emph{For any $m\ge 2$ and $|a|>1$}, \emph{the inequality $\min W(\re(J_{2m-1}(a)[m]))<\min W(\re(J_{2m-1}(a)[2m-1]))$ holds}.

\vspace{3mm}

{\em Proof}. We may assume that $a>1$. Let $A=\re(J_{2m-1}(a)[m])$, $B=\re(J_{2m-1}(a)[2m-1])$, $C=\re J_{m-1}(a)$,
$$D=\frac{1}{2}\left[\begin{array}{cccc} a^{m-1}  & a^m & \cdots & a^{2m-3}\\ a^{m-2} & \ddots & \ddots & \vdots\\
\vdots & \ddots & \ddots & a^m\\ a & \cdots & a^{m-2} & a^{m-1}\end{array}\right],$$
and $b=\min W(B)$. Then
$$A=\left[\begin{array}{cc} C & aD\\ aD^* & C\end{array}\right] \hspace{5mm}\mbox{and}\hspace{5mm} B=\left[\begin{array}{cc} C & D\\ D^* & C\end{array}\right].$$
Since $b=\min W(B)=\min\sigma(B)$ is an eigenvalue of $B$, there is a unit vector $u=u_1\oplus u_2$ in $\mathbb{C}^{m-1}\oplus\mathbb{C}^{m-1}$ such that $Bu=bu$ or
$$\left[\begin{array}{cc} C & D\\ D^* & C\end{array}\right]\left[\begin{array}{c} u_1\\ u_2\end{array}\right]=b\left[\begin{array}{c} u_1\\ u_2\end{array}\right].$$
Hence we have $Cu_1+Du_2=bu_1$ and $D^*u_1+Cu_2=bu_2$. Thus
$$Au=\left[\begin{array}{cc} C & aD\\ aD^* & C\end{array}\right]\left[\begin{array}{c} u_1\\ u_2\end{array}\right]=\left[\begin{array}{c} Cu_1+aDu_2\\ aD^*u_1+Cu_2\end{array}\right]=\left[\begin{array}{c} ((1-a)C+ab)u_1\\ ((1-a)C+ab)u_2\end{array}\right]$$
and, therefore,
\begin{eqnarray}
&&\min W(A)\le\langle Au, u\rangle=\langle ((1-a)C+ab)u_1, u_1\rangle+\langle ((1-a)C+ab)u_2, u_2\rangle \label{eq3}\\
&=& (1-a)(\langle Cu_1, u_1\rangle+\langle Cu_2, u_2\rangle)+ab.\nonumber
\end{eqnarray}
Since $a>1$ and
\begin{eqnarray*}
&& \langle Cu_1, u_1\rangle+\langle Cu_2, u_2\rangle\ge (\min W(C))(\|u_1\|^2+\|u_2\|^2)=\min W(C)\\
&=& \min W(\re J_{m-1}(a))>\min W(\re J_{2m-2}(a))=\min W(\re(J_{2m-1}(a)[2m-1]))\\
&=& \min W(B)=b,
\end{eqnarray*}
where the second inequality follows from Corollary 3.6, we obtain from (\ref{eq3}) that $\min W(A)<(1-a)b+ab=b=\min W(B)$. \hfill $\blacksquare$

\vspace{3mm}

Note that the preceding lemma is not true for $0<|a|\le 1$ as $W(\re (J_3(a)[2]))=W({\scriptsize\left[\begin{array}{cc} 0 & a^2/2\\ \overline{a}^2/2 & 0\end{array}\right]})=[-|a|^2/2, |a|^2/2]$ and $W(\re(J_3(a)[3]))=W({\scriptsize\left[\begin{array}{cc} 0 & a/2\\ \overline{a}/2 & 0\end{array}\right]})=[-|a|/2, |a|/2]$ show.

\vspace{5mm}

{\bf Lemma 3.10.} \emph{For any $m\ge 2$ and $|a|>1$}, \emph{the equality $\min W(\re(J_{2m-1}(a)[m]))=\min W(\re J_{2m-1}(a))$ holds}.

\vspace{5mm}

{\em Proof}. We may assume that $a>1$. Letting $A=\re(J_{2m-1}(a)[m])$, $B=\re J_{2m-1}(a)$ and $b=\min W(A)$, we need only show that $b$ is an eigenvalue of $B$ or, equivalently, $\det(B-bI_{2m-1})=0$. Indeed, if this is the case, then, since $\min\sigma(B[2m-1])$ is between the two smallest eigenvalues of $B$ by the interlacing property of their eigenvalues, we infer from Lemma 3.9 ($b=\min\sigma(A)<\min\sigma(B[2m-1])$) and the above ($b\in\sigma(B)$) that $b=\min\sigma(B)=\min W(B)$.

\vspace{5mm}

To prove $\det(B-bI_{2m-1})=0$, let
$$C=\frac{1}{2}\left[\begin{array}{cccc} a^{m}  & a^{m+1} & \cdots & a^{2m-2}\\ a^{m-1} & \ddots & \ddots & \vdots\\
\vdots & \ddots & \ddots & a^{m+1}\\ a^2 & \cdots & a^{m-1} & a^{m}\end{array}\right].$$
We first transform the rank-1 $(m-1)$-by-$(m-1)$ matrix $C$ into a simpler form via a unitary equivalence. If $U_1$ is the $(n-1)$-by-$(n-1)$ unitary matrix
$$\left[\begin{array}{cccc} 0  &   &     & -1\\   &    & -1 & \\ &  \cdot^{\displaystyle \cdot^{\displaystyle \cdot}} & & \\ -1 &   &   & 0\end{array}\right],$$
then $CU_1=(-1/2)xx^T$, where $x=[a^{m-1} \ a^{m-2} \ \ldots \ a]^T$. Let $c=\sum_{j=1}^{m-1}a^{2j}$, and let $x_1\equiv(1/\sqrt{c})x, x_2, \ldots, x_{m-1}$ be an orthonormal basis of $\mathbb{C}^{m-1}$. If $U_2$ denotes the $(m-1)$-by-$(m-1)$ matrix $[x_1 \ x_2 \ \ldots \ x_{m-1}]$, then $U_2$ is unitary and $U_2^*(CU_1)U_2=\dia(-c/2, 0, \ldots, 0)$. Let $U=U_2\oplus[1]\oplus U_1U_2$, a unitary matrix of size $2m-1$, we have
\begin{equation}
U^*(B-bI_{2m-1})U=\left[\begin{array}{c|c|c} D & \begin{array}{c} \sqrt{c}/2 \\ 0\\ \vdots\\ 0\end{array} & \begin{array}{cccc}-c/2 & & & \\ & 0 & & \\ & & \ddots & \\ & & & 0\end{array}  \\ \hline\begin{array}{cccc} \sqrt{c}/2 & 0 & \cdots & 0\end{array} & -b & \begin{array}{cccc} -\sqrt{c}/2 & 0 & \cdots & 0\end{array}\\ \hline\begin{array}{cccc}-c/2 & & & \\ & 0 & & \\ & & \ddots & \\ & & & 0\end{array} & \begin{array}{c} -\sqrt{c}/2 \\ 0\\ \vdots\\ 0\end{array} & D\end{array}\right], \label{eq4}
\end{equation}
where $D=U_2^*(\re(J_{m-1}(a))-bI_{m-1})U_2\equiv[d_{ij}]_{i, j=1}^{m-1}$. Then we delete the $m$th row and $m$th column from both sides of (\ref{eq4}) to obtain
$$\left[\begin{array}{cc} U_2^* & 0\\ 0 & U_2^*U_1^*\end{array}\right](A-bI_{2m-2})\left[\begin{array}{cc} U_2 & 0\\ 0 & U_2U_1\end{array}\right]=\left[\begin{array}{c|c} D & \begin{array}{cccc}-c/2 & & & \\ & 0 & & \\ & & \ddots & \\ & & & 0\end{array}\\ \hline \begin{array}{cccc}-c/2 & & & \\ & 0 & & \\ & & \ddots & \\ & & & 0\end{array} & D\end{array}\right]\equiv E.$$
This shows the unitary similarity of $A-bI_{2m-2}$ and $E$. Since $b=\min W(A)=\min\sigma(A)$ is in $\sigma(A)$, we have $\det E=\det(A-bI_{2m-2})=0$. Also, note that
$$b=\min W(A)<\min W(B[2m-1])<\min W(\re J_{m-1}(a))=\min\sigma(\re J_{m-1}(a)),$$
where the first inequality is by Lemma 3.9 and the second by Corollary 3.6. This yields the positive-definiteness of $D$. In particular, $D$ is invertible. Hence we can apply the Schur decomposition of $E$ to obtain
\begin{eqnarray}
&& 0=\det E \nonumber\\[3mm]
&=& \det D\cdot\det\Bigg(D-\left[\begin{array}{cccc} -c/2 & & & \\ & 0 & & \\ & & \ddots & \\ & & &  0\end{array}\right]D^{-1}\left[\begin{array}{cccc}-c/2 & & & \\ & 0 & & \\ & & \ddots & \\ & & & 0\end{array}\right]\Bigg) \nonumber\\[3mm]
&=& \det D\cdot\det\left[\begin{array}{cccc} d_{11}-(c^2/4)d & d_{12} & \cdots & d_{1, m-1}\\ d_{21} & d_{22}  & \cdots & d_{2, m-1}\\ \vdots & \vdots &   & \vdots \\ d_{m-1, 1} & d_{m-1, 2} & \cdots & d_{m-1, m-1}\end{array}\right], \label{eq5}
\end{eqnarray}
where $d$ is the $(1, 1)$-entry of $D^{-1}$. It follows that the second determinant of (\ref{eq5}) equals 0. Since this determinant is also equal to $\det D-(c^2d/4)\det(D[1])$ and $d=\det(D[1])/\det D$, a simple computation yields that $d=2/c$. We now apply the Schur decomposition of $U^*(B-bI_{2m-1})U$ to obtain from (\ref{eq4}) that
\begin{eqnarray*}
&& \det(B-bI_{2m-1})
= \det D\cdot\det\Bigg(\left[\begin{array}{ccccc} -b & -\sqrt{c}/2 & 0 & \cdots & 0 \\   -\sqrt{c}/2 & & & & \\ 0 & &D & &\\ \vdots & & & &\\ 0 & & & &\end{array}\right]-\\[3mm]
&& \hspace{3mm}\left[\begin{array}{cc} \sqrt{c}/2 &  0 \ \cdots \ 0 \\  \begin{array}{c} -c/2 \\ 0 \\ \vdots \\ 0\end{array} &   0\end{array}\right]D^{-1}\left[\begin{array}{cc} \sqrt{c}/2 &  -c/2 \ 0 \ \cdots \ 0 \\  \begin{array}{c}  0 \\ \vdots \\ 0\end{array} &  0\end{array}\right]\Bigg)\end{eqnarray*}
\begin{eqnarray*}
\hspace*{7mm} &=& \det D\cdot\det\left[\begin{array}{ccccc} -b-(1/2) &  0 & 0 & \cdots & 0 \\  0& d_{11}-(c/2) & d_{12} & \cdots & d_{1, m-1}\\ 0 &d_{21} & d_{22}  & \cdots & d_{2, m-1}\\ \vdots & \vdots & \vdots &   & \vdots \\ 0 & d_{m-1, 1} & d_{m-1, 2} & \cdots & d_{m-1, m-1}\end{array}\right]\\
&=& -\det D\cdot(b+\frac{1}{2})\cdot\det\left[\begin{array}{cccc} d_{11}-(1/d) & d_{12} & \cdots & d_{1, m-1}\\ d_{21} & d_{22}  & \cdots & d_{2, m-1}\\ \vdots & \vdots &   & \vdots \\ d_{m-1, 1} & d_{m-1, 2} & \cdots & d_{m-1, m-1}\end{array}\right]=0. \hspace{8mm} \blacksquare
\end{eqnarray*}

\vspace{3mm}

{\bf Lemma 3.11.} \emph{For $n\ge 2$ and $a$ in $\mathbb{C}$}, \emph{$-1/2$ is an eigenvalue of $\re J_n(a)$ if and only if $|a|=1$}.

\vspace{3mm}

{\em Proof}. Letting $A_n=2(\re J_n(a))+I_n$, we will prove $\det A_n=(1-|a|^2)^{n-1}$ by induction. Our assertion then follows immediately.

\vspace{3mm}

The asserted expression for $\det A_n$ is obviously true for $n=2$. In general, we multiply the second column of
$$A_n=\left[\begin{array}{cccc} 1  & a  & \cdots & a^{n-1}\\ \overline{a} & 1 & \ddots & \vdots\\
\vdots & \ddots & \ddots & a \\ \overline{a}^{n-1} & \cdots & \overline{a} & 1\end{array}\right]$$
by $-\overline{a}$ and add it to the first column to obtain the matrix
$$B_n=\left[\begin{array}{ccccc} 1-|a|^2  & a  & a^2 & \cdots & a^{n-1}\\ 0 & 1 & a & \cdots & a^{n-2}\\ 0 & \overline{a} & 1 & \ddots & \vdots\\ \vdots & \vdots & \ddots & \ddots & a \\ 0 & \overline{a}^{n-2} & \cdots & \overline{a} & 1\end{array}\right].$$
Hence $\det A_n=\det B_n=(1-|a|^2)\det A_{n-1}$. Then $\det A_n=(1-|a|^2)^{n-1}$ follows by induction. \hfill $\blacksquare$

\vspace{3mm}

{\bf Lemma 3.12.} \emph{For any $n\ge 3$}, $|a|\neq 0, 1$, \emph{and real} $\theta$, \emph{let $b=\max\sigma(\re(e^{i\theta}J_n(a)))$ and let $x=[x_1 \ \ldots x_n]^T$ be a unit vector in $\mathbb{C}^n$ such that $\re(e^{i\theta}J_n(a))x=bx$}. \emph{If $x_j=0$ for some $j$}, $2\le j\le n-1$, \emph{then}

(a) $e^{i\theta}=-1$, \emph{and}

(b) $x_{j-k}=-x_{j+k}$ \emph{for all} $k$, $1\le k\le\min\{j-1, n-j\}$.

\vspace{3mm}

{\em Proof}. We may assume that $0<a\neq 1$.

\vspace{3mm}

(a) Let $A=\re(e^{i\theta}J_n(a))$, $u=[x_1 \ \ldots \ x_{j-1}]^T$, $v=[x_{j+1} \ \ldots \ x_n]^T$, and
$$A[j]=\left[\begin{array}{cc} B & (e^{i\theta}/2)C\\  (e^{-i\theta}/2)C^* & D\end{array}\right],$$
where $B=\re(e^{i\theta}J_{j-1}(a))$, $C=[a^{j-1} \ \ldots \ a^2 \ a]^T[a \ a^2 \ \ldots \ a^{n-j}]$, and $D=\re(e^{i\theta}J_{n-j}(a))$. >From $Ax=bx$, we obtain $A[j]{\scriptsize\left[\begin{array}{c} u \\ v\end{array}\right]}=b{\scriptsize\left[\begin{array}{c} u \\ v\end{array}\right]}$ and, in particular,
\begin{equation}
Bu+\frac{e^{i\theta}}{2}Cv=bu. \label{eq6}
\end{equation}
On the other hand, the equality of the $j$th components of the vectors $Ax$ and $bx$ yields
\begin{equation}
e^{-i\theta}\sum_{\ell=1}^{j-1}a^{j-\ell}x_{\ell}=-e^{i\theta}\sum_{\ell=j+1}^na^{\ell-j}x_{\ell}\equiv c. \label{eq7}
\end{equation}
Therefore,
\begin{equation}
\frac{e^{i\theta}}{2}Cv=\frac{e^{i\theta}}{2}\left[\begin{array}{c} a^{j-1}\\ \vdots \\ a^2\\ a\end{array}\right][ a \ a^2 \ \ldots \ a^{n-j}]v=-\left[\begin{array}{c} a^{j-1}\\ \vdots \\ a^2\\ a\end{array}\right]\frac{c}{2}. \label{eq8}
\end{equation}
Combining (\ref{eq6}), (\ref{eq7}) and (\ref{eq8}) yields
\begin{equation}
\langle (bI_{j-1}-B)u, u\rangle=\frac{e^{i\theta}}{2}\langle Cv, u\rangle=-\frac{c}{2}\langle \left[\begin{array}{c} a^{j-1}\\ \vdots \\ a^2\\ a\end{array}\right], u\rangle=-\frac{c}{2}\cdot \overline{ce^{i\theta}}=-\frac{1}{2}e^{-i\theta}|c|^2.  \label{eq9}
\end{equation}
Since $W(J_{j-1}(a))$ is contained in the interior of $W(J_n(a))$ by Corollary 3.6, we have $\max\sigma(B)=\max W(B)<b$ and hence $bI_{j-1}-B$ is positive-definite. If $u=0$, then, instead of (\ref{eq6}), we obtain $Dv=bv$ from $A[j]{\scriptsize\left[\begin{array}{c} 0 \\ v\end{array}\right]}=b{\scriptsize\left[\begin{array}{c} 0 \\ v\end{array}\right]}$. As above, we may derive that $bI_{n-j}-D$ is positive-definite. It then follows that $v=0$, which, together with $u=0$, gives $x=0$, contradicting our assumption. Hence we must have $u\neq 0$. Then the strict positivity of (\ref{eq9}) yields that $e^{i\theta}=-1$.

\vspace{3mm}

(b) To prove our assertion, note that the equality of the $(j-k)$th components of $Ax$ and $bx$ yields that
\begin{equation}
-\frac{1}{2a^k}\sum_{\ell=1}^{j-k-1}a^{j-\ell}x_{\ell}-\frac{1}{2}\sum_{\ell=j-k+1}^{j-1}a^{\ell-j+k}x_{\ell}
-\frac{a^k}{2}\sum_{\ell=j+1}^na^{\ell-j}x_{\ell}=bx_{j-k}. \label{eq10}
\end{equation}
Similarly, from the $(j+k)$th components of $Ax$ and $bx$, we obtain
\begin{equation}
-\frac{a^k}{2}\sum_{\ell=1}^{j-1}a^{j-\ell}x_{\ell}-\frac{1}{2}\sum_{\ell=j+1}^{j+k-1}a^{j+k-\ell}x_{\ell}
-\frac{1}{2a^k}\sum_{\ell=j+k+1}^na^{\ell-j}x_{\ell}=bx_{j+k}. \label{eq11}
\end{equation}
If $k=1$, then, using (\ref{eq7}) (with $e^{i\theta}=-1$), we can simplify (\ref{eq10}) and (\ref{eq11}) as
\begin{equation}
\frac{1}{2a}(c+ax_{j-1})-\frac{a}{2}c=bx_{j-1} \label{eq12}
\end{equation}
and
\begin{equation}
\frac{a}{2}c-\frac{1}{2a}(c-ax_{j+1})=bx_{j+1}, \label{eq13}
\end{equation}
respectively. Adding (\ref{eq12}) and (\ref{eq13}) and simplifying the resulting equality, we obtain $(b-(1/2))(x_{j-1}+x_{j+1})=0$. Since $b\neq 1/2$ by Lemma 3.11, we have $x_{j-1}=-x_{j+1}$ as asserted.

\vspace{3mm}

We next assume that $x_{j-k}=-x_{j+k}$ for $1\le k<k_0$ and proceed to prove $x_{j-k_0}=-x_{j+k_0}$. Adding (\ref{eq10}) and (\ref{eq11}) results in
\begin{align}
\mbox{  \hspace{8mm}    }&-\frac{1}{2a^{k_0}}\left(\sum_{\ell=1}^{j-k_0}a^{j-\ell}x_{\ell}+\sum_{\ell=j+k_0}^na^{\ell-j}x_{\ell}\right)
-\frac{1}{2}\left(\sum_{\ell=j-k_0+1}^{j-1}a^{\ell-j+k_0}x_{\ell}+\sum_{\ell=j+1}^{j+k_0-1}a^{j+k_0-\ell}x_{\ell}\right) \notag\\[-4mm]
& \label{eq14}\\[-4mm]
&-\frac{a^k_0}{2}\left(\sum_{\ell=j+1}^na^{\ell-j}x_{\ell}+\sum_{\ell=1}^{j-1}a^{j-\ell}x_{\ell}\right)=(b-\frac{1}{2})(x_{j-k_0}+x_{j+k_0}).
\notag\end{align}
Since $x_{j-k}=-x_{j+k}$ for $1\le k\le k_0-1$, the expression within the first parentheses on the left-hand side of (\ref{eq14}) is the same as $\sum_{\ell=1}^{j-1}a^{j-\ell}x_{\ell}+\sum_{\ell=j+1}^na^{\ell-j}x_{\ell}$, which is equal to 0 by (\ref{eq7}). The second term on the left of (\ref{eq14}) is 0 by our assumption that $x_{j-k}=-x_{j+k}$ for $1\le k\le k_0-1$, and the third term is 0 by (\ref{eq7}). Therefore, the right-hand side of (\ref{eq14}) is also 0. Since $b\neq 1/2$ by Lemma 3.11, we conclude that $x_{j-k_0}=-x_{j+k_0}$. \hfill $\blacksquare$

\vspace{3mm}

We are now ready to prove Theorem 3.8.

\vspace{3mm}

{\em Proof of Theorem $3.8$}. Let $A=J_n(a)$. If $n$ is odd, $j=(n+1)/2$ and $|a|>1$, then $b$ is in $\partial W(A)\cap\partial W(A[j])$ by Lemma 3.10. That $b$ is the only element in this intersection is a consequence of Lemma 3.12 (a). Indeed, if $b'$ is any other point in the intersection, then let $\theta$ in $\mathbb{R}$ be such that $e^{i\theta}b'=\max W(\re(e^{i\theta}A))=\max W(\re(e^{i\theta}A[j]))$. Let $x'=[x_1' \ \ldots \ x_{n-1}']^T$ be a unit vector in $\mathbb{C}^{n-1}$ such that $e^{i\theta}b'=\langle (\re(e^{i\theta}A[j]))x', x'\rangle$, and let $x=[x_1' \ \ldots \ x_{j-1}' \ \tb{0}{j\mbox{th}} \ x_{j}' \ \ldots \ x_{n-1}']^T$. Then $x$ is a unit vector in $\mathbb{C}^n$ with $\langle (\re(e^{i\theta} A))x, x\rangle=e^{i\theta}b'$. It then follows that $(\re(e^{i\theta}A))x=e^{i\theta}b'x$. Thus Lemma 3.12 (a) implies that $e^{i\theta}=-1$ and hence $b'=\min W(\re A)=\min\sigma(\re A)=b$. The existence of the unit vector $x$ in $\mathbb{C}^n$ in the asserted form satisfying $(\re J_n(a))x=bx$ is the consequence of the above and Lemma 3.12. We then have $\re\langle J_n(a)x, x\rangle=\langle (\re J_n(a))x, x\rangle=b$. Since, under $|a|>1$, $\partial W(J_n(a))$ has no (vertical) line segment by Proposition 2.11, we infer that $\langle J_n(a)x, x\rangle=b$.

\vspace{3mm}

For the remaining cases, if $0<|a|<1$, then consider $A=f(J_n(a))$, where $f(z)=((1-|a|^2)/a)z-\overline{a}$. By Lemma 2.4, $A$ is of class $S_n$ and hence, by \cite[Theorem 3.3 (b)]{11}, $\partial W(A)\cap\partial W(A[j])=\emptyset$ for all $j$, $1\le j\le n$. It follows easily that $\partial W(J_n(a))\cap\partial W(J_n(a)[j])=\emptyset$ for all $j$. If $n=2$, then our assertion is obviously true. In the following, we assume that $n\ge 3$, $j\neq(n+1)/2$, $|a|>1$ and $\partial W(A)\cap\partial W(A[j])\neq \emptyset$. By considering $U^*AU$ instead of $A$ if necessary, where $U$ is the $n$-by-$n$ unitary matrix
$$\left[\begin{array}{ccc} 0 & & 1\\ &\cdot^{\displaystyle \cdot^{\displaystyle \cdot}} & \\1 & & 0\end{array}\right],$$
we may assume that $j<(n+1)/2$. If $b$ is in $\partial W(A)\cap\partial W(A[j])$, then, as shown in the preceding paragraph, we obtain, using Lemma 3.12, that there is a unit vector $x=[x_1 \ \ldots \ x_{j-1} \ 0 \ -x_{j-1} \ \ldots \ -x_1 \ x_{2j} \ \ldots \ x_n]^T$ in $\mathbb{C}^n$ such that $(\re A)x=bx$. Let $\re A={\scriptsize\left[\begin{array}{cc} B & C\\ C^* & D\end{array}\right]}$, where $B=\re J_{2j-1}(a)$,
$$C=\frac{1}{2}\left[\begin{array}{cccc} a^{2j-1} & a^{2j} & \cdots & a^{n-1}\\ a^{2j-2} & a^{2j-1} & \cdots & a^{n-2}\\ \vdots & \vdots & & \vdots\\ a & a^2 & \cdots & a^{n-2j+1}\end{array}\right],$$
and $D=\re J_{n-2j+1}(a)$. We have
\begin{equation}
Bu+Cv=bu \label{eq15}
\end{equation}
and
\begin{equation}
C^*u+Dv=bv, \label{eq16}
\end{equation}
where $u=[x_1 \ \ldots \ x_{j-1} \ 0 \ -x_{j-1} \ \ldots \ -x_1]^T$ and $v=[x_{2j} \ \ldots \ x_n]^T$. On the other hand, equating the $j$th components of $(\re A)x$ and $bx$ yields $\sum_{\ell=2j}^na^{\ell-2j}x_{\ell}=0$, which then implies that $Cv=0$. Hence (\ref{eq15}) becomes $Bu=bu$. Since $2j-1<n$, $W(J_{2j-1}(a))$ is contained in the interior of $W(J_n(a))$ by Corollary 3.6. Hence $\min W(\re J_{2j-1}(a))>\min W(\re J_n(a))$ or $B>bI_{2j-1}$. It follows from above that $u=0$. Then (\ref{eq16}) becomes $Dv=bv$. A similar argument as above shows that $D>bI_{n-2j+1}$, which yields $v=0$. This means that $x=0$, contradicting our assumption. We conclude that, in this case, we also have $\partial W(A)\cap\partial W(A[j])=\emptyset$. \hfill $\blacksquare$

\vspace{3mm}

{\bf Corollary 3.13.} \emph{Let $A$ be an $S_n^{-1}$-matrix} ($n\ge 2$) \emph{represented as the upper-triangular form $[a_{ij}]_{i, j=1}^n$}, \emph{where $a_{ii}=\lambda$ for all $i$ and $a_{ij}=\overline{\lambda}^{j-i-1}(|\lambda|^2-1)$ if $i<j$}, \emph{and $0$ if $i>j$}, \emph{and let $1\le j_0\le n$}. \emph{Then $\partial W(A)\cap\partial W(A[j_0])\neq \emptyset$ if and only if $n$ is odd and $j_0=(n+1)/2$}. \emph{In this case}, \emph{the intersection consists of a single point $b$ for which there is a unit vector $x=[x_1 \ \ldots \ x_{j_0-1} \ 0 \ -x_{j_0-1} \ \ldots \ -x_1]^T$ in $\mathbb{C}^n$ such that $\langle Ax, x\rangle=b$}.

\vspace{3mm}

{\em Proof}. Since $(\overline{\lambda}/(|\lambda|^2-1))(A-\lambda I_n)=J_n(\overline{\lambda})$ by Lemma 2.4, our assertions follow from Theorem 3.8 and Lemma 3.12. \hfill $\blacksquare$

\vspace{3mm}

We end this paper with a final remark. From \cite[Lemma 2.6]{4}, we can deduce that if $A$ is any $S_n^{-1}$-matrix represented in the standard upper-triangular form as in \cite[Theorem 2.4]{4}, then $\partial W(A)\cap\partial W(A[1])=\partial W(A)\cap\partial W(A[n])= \emptyset$. On the other hand, it can be shown that the $S_3^{-1}$-matrix
$$A=\left[\begin{array}{ccc} 2 & 2\sqrt{3} & 6-12i\\ 0 & 1+2i & 4\sqrt{3}\\ 0 & 0 & 2-3i\end{array}\right]$$
satisfies $\partial W(A)\cap\partial W(A[j])=\emptyset$ for $j=1, 2, 3$. Thus the preceding corollary no longer holds for general $S_n^{-1}$-matrices.

\vspace{1cm}

\noindent
H.-L. Gau, Dept. of Mathematics, National Central University, Chungli 32001, Tai-\\ \hspace*{10mm}wan;  \emph{e-mail}: hlgau@math.ncu.edu.tw

\noindent
P. Y. Wu, Dept. of Applied Mathematics, National Chiao Tung University, Hsinchu\\ \hspace*{10mm}30010, Taiwan; \emph{e-mail}: pywu@math.nctu.edu.tw


\newpage


\begin{thebibliography}{99}

\bibitem{1} H. Bercovici, \textit{Operator Theory and Arithmetic in $H^{\infty}$}, Amer. Math. Soc., Providence, 1988.

\bibitem{2} H. Gaaya, On the numerical radius of the truncated adjoint shift, \textit{Extracta Math.}, 25 (2010), 165--182.

\bibitem{3}  H. Gaaya, A sharpened Schwarz--Pick operatorial inequality for nilpotent operators, arXiv: 1202.3962v1.

\bibitem{4} H.-L. Gau, Numerical ranges of reducible companion matrices, \textit{Linear Algebra Appl.}, 432 (2010), 1310--1321.

\bibitem{5} H.-L. Gau and P. Y. Wu, Numerical range of $S(\phi)$, \textit{Linear Multilinear Algebra}, 45 (1998), 49--73.

\bibitem{6} H.-L. Gau and P. Y. Wu, Dilation to inflations of $S(\phi)$, \textit{Linear Multilinear Algebra}, 45 (1998), 109--123.

\bibitem{7} H.-L. Gau and P. Y. Wu, Lucas' theorem refined, \textit{Linear Multilinear Algebra}, 45 (1999), 359--373.

\bibitem{8} H.-L. Gau and P. Y. Wu, Conditions for the numerical range to contain an elliptic disc, \textit{Linear Algebra Appl.}, 364 (2003), 213--222.

\bibitem{9} H.-L. Gau and P. Y. Wu, Companion matrices: reducibility, numerical ranges and similarity to contractions, \textit{Linear Algebra Appl.}, 383 (2004), 127--142.

\bibitem{10} H.-L. Gau and P. Y. Wu, Numerical ranges of nilpotent operators, \textit{Linear Algebra Appl.}, 429 (2008), 716--726.

\bibitem{11} H.-L. Gau and P. Y. Wu, Numerical ranges and compressions of $S_n$-matrices, \textit{Oper. Matrices}, to appear.

\bibitem{12} U. Grenander and G. Szeg\H{o}, \textit{Toeplitz Forms and Their Applications}, Univ. of California Press, Berkeley, 1958.

\bibitem{13} U. Haagerup and P. de la Harpe, The numerical radius of a nilpotent operator on a Hilbert space, \textit{Proc. Amer. Math. Soc.}, 115 (1992), 371--379.

\bibitem{14} R. A. Horn and C. R. Johnson, \textit{Topics in Matrix Analysis}, Cambridge Univ. Press, Cambridge, 1991.

\bibitem{15} M. Kac, W. L. Murdock and G. Szeg\H{o}, On the eigenvalues of certain Hermitian forms, \textit{J. Rational Mech. Anal.}, 2 (1953), 767--800.

\bibitem{16} D. S. Keeler, L. Rodman and I. M. Spitkovsky, The numerical range of $3\times 3$ matrices, \textit{Linear Algebra Appl.}, 252 (1997), 115--139.

\bibitem{17} R. Kippenhahn, \"{U}ber den Wertevorrat einer Matrix,
\textit{Math. Nachr.}, 6 (1951), 193--228. (English translation: P. F. Zachlin
and M. E. Hochstenbach, On the numerical range of a matrix, \textit{Linear
Multilinear Algebra}, 56 (2008), 185--225.)

\bibitem{18} F. Kirwan, \textit{Complex Algebraic Curves}, Cambridge Univ. Press, Cambridge, 1992.

\bibitem{19} C.-K. Li, B.-S. Tam and P. Y. Wu, The numerical range of a nonnegative matrix, \textit{Linear Algebra Appl.}, 350 (2002), 1--23.

\bibitem{20} J. Maroulas, P. J. Psarrakos and M. J. Tsatsomeros, Perron--Frobenius type results on the numerical range, \textit{Linear Algebra Appl.}, 348 (2002), 49--62.

\bibitem{21} D. Sarason, Generalized interpolation in $H^{\infty}$, \textit{Trans. Amer. Math. Soc.}, 127 (1967), 179--203.

\bibitem{22} B. Sz.-Nagy, C. Foia\c{s}, H. Bercovici and L. K\'{e}rchy,  \textit{Harmonic Analysis of Operators on Hilbert Space}, 2nd ed., Springer, New York, 2010.

\bibitem{23} B.-S. Tam and S. Yang, On matrices whose numerical ranges have circular or weak circular symmetry, \textit{Linear Algebra Appl.}, 302/303 (1999), 193--221.

\bibitem{24} K.-Z. Wang and P. Y. Wu, Diagonals and numerical ranges of weighted shift matrices, \textit{Linear Algebra Appl.}, (2012), http://dx.doi.org/10.1016/j.laa.2012.08.007.

\end{thebibliography}
\end{document}